\documentclass{amsart}
\usepackage[utf8]{inputenc}
\usepackage{xfrac, url, longtable, placeins, mdwtab, amsfonts, amsmath, todonotes}
\usepackage{lscape}
\newtheorem{thm}{Theorem}[section]

\newtheorem{prop}[thm]{Proposition}
\theoremstyle{definition}
\newtheorem{rem}[thm]{Remark}
\newtheorem{alg}{Algorithm}
\newcommand{\gd}{\mathfrak{d}}

\newcommand{\gp}{\mathfrak{p}}

\newcommand{\CL}{\mathcal{L}}
\newcommand{\CO}{\mathcal{O}}

\newcommand{\CS}{\mathcal{S}}
\newcommand{\FF}{\mathbb{F}}
\newcommand{\NN}{\mathbb{N}}
\newcommand{\QQ}{\mathbb{Q}}

\newcommand{\ZZ}{\mathbb{Z}}

\newcommand{\QQt}{\QQ(\vartheta)}
\newcommand{\form}[1]{\langle #1\rangle}
\newcommand{\st}{\bigm|}
\newcommand{\nilrad}{\mathfrak{N}}

\newenvironment{poc}{\begin{proof}[Proof of correctness]}{\end{proof}}
\newcommand{\term}[1]{\emph{#1}}
\newcommand{\card}{\sharp}
\newcommand{\disc}{\mbox{disc}\,}
\newcommand{\ind}{\mbox{ind}\,}
\newcommand{\ord}{\mbox{ord}\,}
\newcommand{\sgn}{\mbox{sgn}\,}

\begin{document}
\title{Computing with quadratic forms over number fields}
\author{Przemys{\l}aw Koprowski}
\address{Faculty of Mathematics\\
  University of Silesia\\
  ul. Bankowa 14\\ 
  PL-40-007 Katowice, Poland}
\email{pkoprowski@member.ams.org}
\author{Alfred Czoga{\l}a}
\address{Faculty of Mathematics\\
  University of Silesia\\
  ul. Bankowa 14\\ 
  PL-40-007 Katowice, Poland}
\email{alfred.czogala@us.edu.pl}
\maketitle

\begin{abstract}
This paper presents fundamental algorithms for the computational theory of quadratic forms over number fields. In the first part of the paper, we present algorithms for checking if a given non-degenerate quadratic form over a fixed number field is either isotropic (respectively locally isotropic) or hyperbolic (respectively locally hyperbolic). Next we give a method of computing the dimension of an anisotropic part of a quadratic forms. The second part of the paper is devoted to algorithms computing two field invariants: the level and the Pythagoras number. Ultimately we present an algorithm verifying whether two number fields have isomorphic Witt rings (i.e. are Witt equivalent).
\end{abstract}


\section{Introduction}
The algebraic theory of quadratic forms is a mature and important branch of mathematics. Yet still, the computational side of this theory is seriously under-developed. The majority of research concentrate on forms over the rationals. Consequently, while over~$\QQ$ there already a couple of algorithms for solving a highly non-trivial problem of determining isotropic vectors of a quadratic forms (see e.g. \cite{CreRus03, Simon05, Castel13}), little has been done so far for forms over number fields (i.e. finite extensions of $\QQ$). The algebraic theory of quadratic forms over number fields are very like the theory over the rationals, nevertheless the computational approach seems to be rudimentary here. The aim of this article is to partially fill this gap, as well as provoke further discussion and future research.

This paper is organized as follows: in Section~\ref{sec_isotropy} we present an algorithm (see Algorithm~\ref{alg_iso_global}) checking if a given form (over a fixed number field $K$) is isotropic. This algorithm uses sub-procedures (Algorithms~\ref{alg_iso_nd} and~\ref{alg_iso_d}) deciding whether the form is isotropic at a non-archimedean prime of $K$ (respectively odd or even). These two algorithm may be of an independent interest to the reader. Next, in Section~\ref{sec_hiperbolic} we show Algorithm~\ref{alg_hyp_global} determining if a quadratic form is hyperbolic, again utilizing the local approach. 

It is known that any non-degenerate form can be uniquely decomposed into an orthogonal sum of its anisotropic part and a hyperbolic form (one of these two parts may of course be void if the form in question is either anisotropic or hyperbolic itself). In Section~\ref{sec_Witt_index} we shows a procedure computing the dimension of an anisotropic part of a quadratic form.

In Sections~\ref{sec_level}--\ref{sec_Witt_eq} we go a step further and develop algorithms computing invariants of the ground fields, that play important roles in the algebraic theory of quadratic forms. Algorithm~\ref{alg_level} computes the level $s(K)$ of a number field $K$, which is the length of the shortest representation of $-1$ as a sum of squares. Another invariant of the field is the minimal number of squares needed to represent any sum of squares. This invariant is called the Pythagoras number and is computed by Algorithm~\ref{alg_Pythagoras_number}.

Recall that the set $WK$ of similarity classes of non-degenerate symmetric bilinear forms over a given base field $K$ is a ring with operations induced by the orthogonal sum and the tensor product. It is called the \emph{Witt ring} of the field $K$. Because a bilinear form defines an orthogonal geometry on the vector space on which it is defined, thus the Witt ring can be viewed as an algebraic structure encoding information on all possible orthogonal geometries over a given base field. Two fields are said to be \emph{Witt equivalent}, if their Witt rings are isomorphic. The set of global field invariants that fully determine its Witt equivalence class was described in \cite{Szymiczek91}. In Section~\ref{sec_Witt_eq} we present Algorithm~\ref{alg_WE_invariants} computing all these invariants. In particular the algorithm may be used to verify whether two number fields are Witt equivalent.

The authors implemented all the algorithms presented in this paper in a computer algebra system Sage\nocite{sage}. Using this implementation, we were able to find representatives of Witt classes of number fields of low degrees. These results are presented in Tables~\ref{tbl_cubic}--\ref{tbl_sextic}. Moreover, using our implementation, we were able to give an affirmative answer to Conner's question for number field of degree not exceeding $6$ (for details see the last section of the paper). 

In these paper, $K=\QQt$ is always a number field specified by the minimal polynomial of $\vartheta$ over $\QQ$ and $\CO_K$ is the integral closure of $\ZZ$ in $K$. Two basic building blocks that we use in subsequent algorithms are procedures that test whether a given algebraic number $a\in K$ is a square: either in its base field $K$ or in a completion $K_\gp$, where $\gp$ is a prime of $K$. A procedure testing whether an element is a square in a number field is available as standard in computer algebra systems. On the other hand, testing whether $a$ is a square in a completion~$K_\gp$ is obviously equivalent to testing whether $x^2-a$ is irreducible in $K_\gp[x]$. There are known algorithms for testing irreducibility of a polynomial in local fields. These include for example: Montes' algorithm (see e.g. \cite{Veres09} or \cite{GMN11, GMN12}) or variations of Zassenhaus Round Four algorithm (see e.g. \cite{Pauli01, Pauli10}).

In the algorithms presented below, an input is a non-degenerate diagonal quadratic form with coefficients in some number field $K$. Since $K$ is the field of fractions of $\CO_K$ and for every $a,b\in \CO_K$, both $\sfrac{a}{b}$ and $a\cdot b$ belong to the same square-class on $\sfrac{\dot{K}}{\dot{K}^2}$, hence in Algorithms \ref{alg_iso_ff}--\ref{alg_aniso_part_dim} we usually assume that the coefficients of the quadratic form come from~$\CO_K$.

\section{Isotropy of a quadratic form}\label{sec_isotropy}
In this section, we present an algorithm that checks if a given form $\varphi$ over a number field $K$ is isotropic or not. The organization of this section reflects the general idea of solving the problem locally. Hence, Algorithms~\ref{alg_iso_nd}, \ref{alg_iso_d} and~\ref{alg_signatures} deal respectively with odd and even finite fields and real infinite primes of $K$. Finally, Algorithm~\ref{alg_iso_global} checks if the form is globally isotropic, using the above-mentioned algorithms as sub-procedures. 

Below we utilize the notion of the discriminant of a quadratic form. Recall (see e.g. \cite[Definition~15.2.1]{Szymiczek91}) that for a quadratic form $\varphi$, we define the discriminant of $\varphi$ by the formula
\[
\disc\varphi := (-1)^{\sfrac{d(d-1)}{2}} \det\varphi,
\]
where $d= \dim \varphi$.

\begin{alg}\label{alg_iso_ff}
Let $\gp$ be an odd prime of $K$ and $\varphi = \form{a_1, \dotsc, a_d}$ be a non-degenerate diagonal quadratic form with all its entries being $\gp$-adic units. This algorithm returns true if and only if the residual form $\varphi \otimes K/\gp$ is isotropic, otherwise it returns false.
\begin{enumerate}
\item If $\dim \varphi = 1$, return false.
\item If $\dim \varphi = 2$, return true when $\disc \varphi$ is a square in $K/\gp$, otherwise return false.
\item If $\dim \varphi>2$, return true.
\end{enumerate}
\end{alg}

The correctness of the above algorithm follows immediately from \cite[Theorem~I.3.2]{Lam05}.

\begin{alg}\label{alg_iso_nd}
Let $\gp$ be an odd prime of a number field $K$. Given a non-degenerate quadratic form $\varphi$, this algorithm returns true if $\varphi\otimes K_\gp$ is isotropic and false otherwise.
\begin{enumerate}
\item If $\dim \varphi = 1$, return false.
\item If $\dim \varphi \geq 5$, return true.
\item\label{st_split} Let $\{a_1, \dotsc, a_d\}$ be the list of coefficients of a diagonalization of $\varphi$, all $a_i\in \CO_K$. Partition this list into two sublists depending on the parity of the $\gp$-adic valuation:
\begin{align*}
\varphi_0 &:= \bigl\{ a_i\cdot \pi^{-\ord_\gp a_i}\st \ord_\gp a_i\equiv 0\pmod{2}\bigr\},\\
\varphi_1 &:= \bigl\{ a_i\cdot \pi^{-\ord_\gp a_i}\st \ord_\gp a_i\equiv 1\pmod{2}\bigr\}.
\end{align*}
Here $\pi$ is a uniformizer of $\gp$ (see Remark~\ref{rem_uniformizer} below).
\item Use Algorithm~\ref{alg_iso_ff} to verify whether any of $\varphi_0$, $\varphi_1$ is isotropic over $K/\gp$. Return true if Algorithm~\ref{alg_iso_ff} returned true at least once, otherwise return false.
\end{enumerate}
\end{alg}

\begin{rem}\label{rem_uniformizer}
In order to find a uniformizer of a given prime in step~\eqref{st_split} of the above algorithm, one may use for example \cite[Algorithm~4.8.17]{Cohen93} or \cite[\S3]{GMN13}.
\end{rem}

The correctness of the algorithm follows from \cite[Proposition~VI.1.9]{Lam05}. Next, we consider even primes. Recall (see e.g. \cite[Definition~V.3.17]{Lam05}) that the Hasse invariant of a quadratic form $\varphi = \form{a_1, \dotsc, a_d}$ at a prime $\gp$ is:
\begin{equation}\label{eq_Hasse}
s_\gp(\varphi) := \prod_{1\leq i<j\leq d} (a_i, a_j)_\gp,
\end{equation}
where $(a_i, a_j)_\gp$ denotes the $\gp$-adic Hilbert symbol. An algorithm for computing the Hilbert symbol in a completion of a number field was recently presented in \cite{Voight13}. We use it to verify whether a quadratic form is isotropic over a dyadic completion of~$K$.

\begin{alg}\label{alg_iso_d}
Let $\gd$ be an even prime of $K$ and $\varphi$ be a non-degenerate quadratic form over $K$. This algorithm returns true if and only if $\varphi\otimes K_\gd$ is isotropic, otherwise it returns false.
\begin{enumerate}
\item\label{step_1} If $\dim \varphi \leq 1$, then return false and quit.
\item\label{step_2} If $\dim\varphi = 2$, then check whether $\disc \varphi$ is a square in~$K_\gd$. If so, then return true and quit, otherwise return false and quit.
\item If $\dim\varphi = 3$, then proceed as follows:
  \begin{enumerate}
  \item Compute the Hilbert symbol $\bigl(-1,-\det(\varphi)\bigr)_\gd$ by applying \cite[Algorithm~6.6]{Voight13}.
  \item Use Eq.~\eqref{eq_Hasse} and \cite[Algorithm~6.6]{Voight13} to compute the Hasse invariant $s_\gd(\varphi)$ of $\varphi$ at $\gd$.
  \item If $\bigl(-1,-\det(\varphi)\bigr)_\gd = s_\gd(\varphi)$, then return true otherwise return false
  \end{enumerate}
\item\label{step_4} If $\dim\varphi = 4$, then proceed as follows:
  \begin{enumerate}
  \item Check if $\det \varphi$ is a square in $K_\gd$. If not, then return true and quit.
  \item\label{step_b} If $\det \varphi\in (K_\gd^\times)^2$, then use Eq.~\eqref{eq_Hasse} and \cite[Algorithm~6.6]{Voight13} to compute the Hasse invariant $s_\gd(\varphi)$ and the Hilbert symbol $(-1,-1)_\gd$. Return true if they are equal, return false if they are not.
  \end{enumerate}
\item\label{step_5} If $\dim\varphi\geq 5$, then return true.
\end{enumerate}
\end{alg}

\begin{poc}
An unary form is never isotropic and a quintic or higher-di\-men\-sio\-nal form over a dyadic field is always isotropic by the means of \cite[Theorem~VI.2.12]{Lam05}. This justifies steps~\eqref{step_1} and \eqref{step_5}. Next, it is well known that a binary form is isotropic if and only if its determinant is a minus square, which proves step \eqref{step_2}. On the other hand, if the form has dimension three, then \cite[Proposition~V.3.22]{Lam05} asserts that it is isotropic if and only if $\bigl(-1,-\det(\varphi)\bigr)_\gd = s_\gd(\varphi)$.

This leaves us with quaternary forms. Now, \cite[Corollary~VI.2.15]{Lam05} asserts that over a local field there is only one anisotropic form of dimension $4$ and its determinant is a square. Thus, if $\det \varphi\notin (K_\gd^\times)^2$, then $\varphi\otimes K_\gd$ is necessarily isotropic. On the other hand, if 
$\det \varphi\in (K_\gd^\times)^2$, then \cite[Proposition~V.3.23]{Lam05} provides us with a needed criterion for isotropy.
\end{poc}

After covering the finite primes we need a tool do deal with the infinite ones, as well. Recall (see e.g. \cite[p.~34]{Lam05}), that the \term{signature} of a non-degenerate quadratic form $\varphi$ with respect to an ordering $\beta$ of the coefficient field is the difference between the number of positive and negative entries of a diagonalization $\form{a_1, \dotsc, a_d}$ of $\varphi$:
\[
\sgn_\beta(\varphi) := \card \{a_i\st a_i >_\beta 0\} - \card \{a_i\st a_i <_\beta 0\}.
\]
This number is known to be independent of a choice of an actual diagonalization of~$\varphi$. We now present an algorithm computing the signatures of the form with respect to all orderings of $K$.

\begin{alg}\label{alg_signatures}
Let $\QQt$ be a number field, specified by a minimal polynomial $f\in \QQ[x]$ of its generator $\vartheta$. Given a non-degenerate diagonal quadratic form $\varphi = \form{a_1, \dotsc, a_d}$ with coefficients in $\CO_K$, this algorithm computes the list of signatures of $\varphi$ with respect to all orderings of~$K$.
\begin{enumerate}
\item Use \cite[Algorithm~10.64]{BPR03} to find a list $(\sigma_1, \dotsc, \sigma_r)$ of Thom encodings of all real roots of $x_1 <\dotsb < x_r$ of the generating polynomial $f$;
\item For every coefficient $a_i$ of $\varphi$ proceed as follows:
  \begin{enumerate}
  \item Let $g_0, \dotsc, g_{n-1}\in \QQ$ be the coordinates of $a_i$ with respect to the power basis $\{ 1, \vartheta, \dotsc, \vartheta^{n-1}\}$ (i.e. $a_i = g(\vartheta)$ with $g = g_0 x+\dotsb + g_{n-1} x^{n-1}\in \QQ[x]$).
  \item Use \cite[Algorithm~10.67]{BPR03} to determine the signs:
  \[ s_{i1} = \sgn g(x_1), \dotsc, s_{ir} = \sgn g(x_r) \]
  of the polynomial $g$ at the roots of $f$.
  \end{enumerate}
 \item Return the list of sums $\bigl(\sum_{i=1}^d s_{i1}, \dotsc, \sum_{i=1}^d s_{ir}\bigr)$.
\end{enumerate}
\end{alg}

The correctness of the algorithm follows immediately from the correctness of \cite[Algorithms~10.64 and~10.67]{BPR03}. 

\begin{rem}
The above algorithm is used subsequently in step~\eqref{step_r} of Algorithms~\ref{alg_iso_global}, \ref{alg_hyp_global} and step~\eqref{st_real_aniso_dim} of Algorithm~\ref{alg_aniso_part_dim}. As an alternative approach, one could use here an interval arithmetic and evaluate the signs of $\bigl(g_j(\vartheta_i)\bigr)$ by the means of \cite[\S8.5, Sign evaluation]{Mishra93}.
\end{rem}

\begin{rem}\label{rem_factor}
In algorithms \ref{alg_iso_global}, \ref{alg_hyp_global}, \ref{alg_aniso_part_dim}, \ref{alg_level} and~\ref{alg_Pythagoras_number} below, we need to perform factorizations of two kinds. The first one is to find all even primes of a given field $K$, i.e. to factor $2\CO_K$. The other one is to find all primes dividing any of the coefficients of a given quadratic form. The factorization of an ideal in a number field corresponds to the factorization of a polynomial in a local field (see a comment at the end of the introduction). Algorithms for the factorization of ideals are known and described in computational algebraic number theory literature. One may refer for example to \cite[\S6.2.5]{Cohen93} and \cite[2.3.22]{Cohen00} or to a newer algorithm described in \cite[\S2.2]{GMN13}.
\end{rem}

Now, we are finally ready to present the main algorithm of this section, that checks if a form is isotropic over a given number field.

\begin{alg}\label{alg_iso_global}
Given a non-degenerate diagonal quadratic form $\varphi = \form{a_1, \dotsc, a_d}$ over $K$ with $a_i\in \CO_K$, this algorithm returns true if and only if $\varphi$ is isotropic and false if it is not.
\begin{enumerate}
\item If $\dim \varphi\leq 1$, then return false and quit.
\item If $\dim\varphi = 2$, then check if $\disc\varphi$ is a square in $K$. If so, then return false; if not, return true.
\item\label{step_r} Use Algorithm~\ref{alg_signatures} to compute the list $(s_1, \dotsc, s_r)$ of the signatures of $\varphi$ under all real embeddings of $K$. If $|s_j|= \dim \varphi$ for any $1\leq j\leq r$, then return false and quit.
\item\label{step_d} Factor $2\CO_K$ into prime ideals $2\CO_K = \gd_1^{e_1}\dotsb \gd_n^{e_n}$ in $\CO_K$ (see Remark~\ref{rem_factor}). For each $\gd_i$ use Algorithm~\ref{alg_iso_d} to check if $\varphi\otimes K_{\gd_i}$ is isotropic. If the algorithm returns false, for at least one $\gd_i$, then return false and quit.
\item\label{step_nd} Find all odd primes $\gp$ of $K$ dividing any of the coefficients $a_i$ of $\varphi$. For each such a prime $\gp$ call Algorithm~\ref{alg_iso_nd}. If the procedure returns false at least once, then return false and quit.
\item Return true.
\end{enumerate}
\end{alg}

\begin{poc}
The cases of unary and binary forms are trivial. For forms of higher dimension we use the local-global principle \cite[Principle VI.3.1]{Lam05}. The form is isotropic over $K$ if and only if it is isotropic over all the completions of $K$. Now $\varphi$, having dimension at least three, is trivially isotropic at all odd primes that do not divide any of the coefficients. These are almost all primes of $K$. Thus, we are left with only finitely many cases to check: finitely many real places treated in step~\eqref{step_r}, finitely many dyadic places covered by step~\eqref{step_d} and finitely many non-dyadic primes considered in step~\eqref{step_nd}.
\end{poc}

\section{Hyperbolicity of a quadratic form}\label{sec_hiperbolic}
In this section we present an algorithm checking another fundamental property of a quadratic form, namely whether it is hyperbolic (hence, a zero element in the Witt group). The general idea is similar to the one adopted in the previous section. Again, we treat the problem locally, separately for finite and real infinite primes of~$K$. 

\begin{alg}\label{alg_hyp_local}
Let $\gp$ be a finite prime of a number field $K$ (either even or odd). Given a non-degenerate quadratic form $\varphi$, this algorithm returns true if the form $\varphi_\gp := \varphi\otimes K_\gp$ is hyperbolic and false otherwise.
\begin{enumerate}
\item If $\dim\varphi$ is odd, then return false and quit.
\item\label{step_local_square} Compute the discriminant $\disc \varphi$ and check if it is a square in the completion $K_\gp$. If it is not a square, then return false and quit.
\item Use Eq.~\eqref{eq_Hasse} and \cite[Algorithm~6.6]{Voight13} to compute the Hasse invariant $s_\gp(\varphi)$ and the power $\smash{(-1,-1)_\gp^{m(m-1)/2}}$ of $\gp$-adic Hilbert symbol, where $2m=\dim \varphi$. Return true if they are equal, return false if they are not.
\end{enumerate}
\end{alg}

\begin{poc}
Take a form $\varphi$ of an even dimension. If the discriminant $\disc \varphi$ is a square in $K_\gp$ and the Hasse invariant $s_\gp(\varphi)$ equals $\smash{(-1,-1)_\gp^{m(m-1)/2}}$, then $\varphi$ is isometric to the hyperbolic space $m\form{1,-1}$ by \cite[Proposition V.3.25]{Lam05}.
\end{poc}

\begin{alg}\label{alg_hyp_global}
Given a non-degenerate diagonal quadratic form $\varphi = \form{a_1, \dotsc, a_d}$ over $K$ with $a_i\in \CO_K$, this algorithm returns true if and only if $\varphi$ is hyperbolic, otherwise it returns false.
\begin{enumerate}
\item If $\dim\varphi$ is odd, then return false and quit.
\item Compute the discriminant $\disc\varphi$. Check if $\disc\varphi$ is a square in $K$. If it is not, then return false.
\item Use Algorithm~\ref{alg_signatures} to compute the list $(s_1, \dotsc, s_r)$ of the signatures of $\varphi$ under all real embeddings of $K$. If $s_j\neq 0$ for any $1\leq j\leq r$, then return false and quit.
\item\label{step_hyp_global_primes} Let $\CL$ be the set consisting of all odd primes of $K$ dividing any of the coefficients $a_i$ of $\varphi$ and of all even primes of~$K$.
\item Apply Algorithm~\ref{alg_hyp_local} to every $\gp\in \CL$ to check if $\varphi\otimes K_\gp$ is hyperbolic. If it returns false, for at least one $\gp$, then return false and quit.
\item Return true.
\end{enumerate}
\end{alg}

\begin{poc}
It is well known that the discriminant of a hyperbolic form is a square and its dimension has to be even.
Moreover, by the well known Weak Hasse Principle, a quadratic form is hyperbolic over a number field if and only if it is hyperbolic over every completion (finite or real infinite) of the field. Over the reals, the form is hyperbolic, when its signature is null. This proves that the algorithm returns true for all hyperbolic forms.

Conversely, suppose that the algorithm returns true for some non-degenerate form $\varphi$. Thus, $\dim\varphi$ is even, its discriminant is a square an it has a zero signature with respect to every ordering of $K$. 

Recall that a quadratic form over a non-dyadic local field $K_\gp$ decomposes into a sum $\varphi \cong \varphi_1\perp \pi\cdot \varphi_2$, where $\pi$ is a $\gp$-adic uniformizer and the coefficients of $\varphi_1, \varphi_2$ are $\gp$-adic units. Recall (see e.g. \cite[\S~VI.1]{Lam05}) that a map $\varphi\mapsto \varphi_2\otimes (\sfrac{K}{\gp})$ is a well defined homomorphism of Witt groups (i.e. additive groups of Witt rings) $WK_\gp\to W(\sfrac{K}{\gp})$ called the \term{second residue homomorphism}.

If the second residual homomorphisms with respect to all odd primes of $K$ are null, then the Witt class of $\varphi$ sits in $\nilrad (W\CO_K) \cap I^2K$ by \cite[Corollary~IV.4.5]{MH1973}, where $\nilrad (W\CO_K)$ denotes the nilradical of the Witt ring of $\CO_K$ and $IK$ is the fundamental ideal of the Witt ring $WK$. Clearly one needs to check only these primes that divide any of the coefficients of $\varphi$ as we do in step~\eqref{step_hyp_global_primes}. Since our algorithm returns true for the form $\varphi$, hence in particular $\varphi\otimes K_\gd$ is hyperbolic, and so $c_\gd(\varphi) = 1$, for every even prime $\gd$ of $K$. Here
\[
c_\gd(\varphi) = (-1, -1)_\gd^{m(m-1)/2} s_\gd(\varphi),\qquad m = \frac{1}{2}\dim\varphi
\]
is the Hasse-Witt invariant of $\varphi$ (c.f. \cite[Proposition~V.3.20]{Lam05}. Now, the map $\varphi\mapsto \bigl(c_{\gd_1}(\varphi), \dotsc, c_{\gd_{g-1}}(\varphi)\bigr)$ is an isomorphism from $\nilrad(W\CO_K)\cap I^2K$ onto $\{\pm1\}^{g-1}$, where $\gd_1, \dotsc, \gd_g$ are all the dyadic primes of $K$, by \cite[Proposition~3.5]{Czogala01}. It follows that the class of $\varphi$ in $WK$ is null, hence $\varphi$ is a hyperbolic form.
\end{poc}

\section{Witt index of a quadratic form}\label{sec_Witt_index}
Recall (see e.g. \cite[Chapter~i, \S4]{Lam05}) that any non-degenerate quadratic form $\varphi$ can be uniquely (up to an isometry) decomposed as $\varphi = \psi \perp H$, where $\psi$ is an anisotropic form, called the \term{anisotropic part} of $\varphi$ and $H$ is hyperbolic. The number of hyperbolic planes constituting $H$ (i.e. half of the dimension of $H$) is called the \term{Witt index} of $\varphi$ and denoted $\ind (\varphi)$. In this section we present an algorithm that computes the dimension of the anisotropic part of $\varphi$. It can be also used to deduce the Witt index since clearly $\ind\varphi = \sfrac12\cdot(\dim\varphi - \dim\psi)$. Again, the problem is first solved locally (see Algorithm~\ref{alg_aniso_part_local_dim}) and then the local solution is used to derive the global one in Algorithm~\ref{alg_aniso_part_dim}.

\begin{alg}\label{alg_aniso_part_local_dim}
Given a non-degenerate quadratic form $\varphi$ over a number field $K$ and a finite prime $\gp$, this algorithm computes the dimension of the anisotropic part of $\varphi_\gp:= \varphi\otimes K_\gp$ over the completion $K_\gp$.
\begin{enumerate}
\item If $\dim\varphi$ is even, proceed as follows:
  \begin{enumerate}
  \item Use Algorithm~\ref{alg_hyp_local} to check if $\varphi_\gp$ is hyperbolic. If it is, then return $0$ and quit.
  \item Check if $\disc \varphi$ is a square in $K_\gp$. If so, then return $2$ and quit. 
  \item Return $4$.
  \end{enumerate}
\item If $\dim\varphi$ is odd, proceed as follows:
  \begin{enumerate}
  \item\label{step_psi} Let $d:= \dim \varphi$ and take $\psi:= \varphi\perp \form{(-1)^{\sfrac{d(d+1)}{2}}\cdot \det \varphi}$.
  \item Use Algorithm~\ref{alg_hyp_local} to check if $\psi\otimes K_\gp$ is hyperbolic. If it is, then return $1$ and quit.
  \item Return $3$.
  \end{enumerate}
\end{enumerate}
\end{alg}

\begin{poc}
First assume that $\varphi$ is an even-dimensional form, so $\varphi_\gp\in IK_\gp$. If it is not hyperbolic, then its class in the Witt ring $WK_\gp$ is not zero. Suppose that $\disc \varphi$ is a square in $K_\gp$. It follows that $\varphi_\gp \in I^2K_\gp$. But for a local field there is only one non-zero element of $I^2K_\gp$, namely the form $\eta_\gp = \form{1, u, \pi, u\pi}$, here $u$ is a $\gp$-adic unit such that $K_\gp(\sqrt{u})$ is the unique unramified extension of $K_\gp$ (see \cite[Corollary~VI.2.15]{Lam05}). It follows that the anisotropic part pf $\varphi_\gp$ has dimension $4$. Conversely, suppose that $\disc \varphi$ is not a square in $K_\gp$. Therefore $\varphi_\gp\in IK_\gp\setminus I^2K_\gp$ and so the anisotropic part of $\varphi_\gp$ has dimension~$2$.

Now assume that the dimension of $\varphi$ is odd. Hence, the form $\psi$ constructed in step~\eqref{step_psi} is an even dimensional form and its discriminant is a square in $K_\gp$. Consequently, the Witt class of $\psi_\gp := \psi\otimes K_\gp$ sits in $I^2K_\gp$. If $\psi_\gp$ is hyperbolic then $\psi_\gp\cong \frac{d+1}{2}\form{1,-1}$, hence $\varphi_\gp\perp \form{1,-1}\cong \form{c}\perp \frac{d+1}{2}\form{1,-1}$ for $c= -(-1)^{\sfrac{d(d+1)}{2}}\det \varphi$. This implies that the anisotropic part of $\varphi_\gp$ is unary. Conversely, suppose that $\psi_\gp$ is not hyperbolic. As in the first part of the proof, this leads to $\psi_\gp = \eta_\gp$ in the Witt ring $WK_\gp$. In particular, the Witt classes of $\varphi_\gp$ and $\form{c,1,u,\pi,u\pi}$ are equal. But a quintic form over a local field is necessarily isotropic and so it is similar to either ternary or unary form. We claim that the unary case is impossible. Indeed, if
\[
\form{c,1,u,\pi,u\pi} \cong \form{x}\perp 2\form{1,-1},
\]
then square classes of $c$ and $x$ are equal and the Witt cancellation theorem asserts that the forms $\form{1,u,\pi,u\pi}$ and $2\form{1,-1}$ are isometric over $K_\gp$ contradicting \cite[Corollary~VI.2.15]{Lam05}. All in all, $\form{c,1,u,\pi,u\pi}$ has a ternary anisotropic part and so has~$\varphi_\gp$.
\end{poc}

\begin{alg}\label{alg_aniso_part_dim}
Given a non-degenerate quadratic form $\varphi$ over a number field $K$, this algorithm computes the dimension of the anisotropic part of $\varphi$.
\begin{enumerate}
\item\label{st_real_aniso_dim} Use Algorithm~\ref{alg_signatures} to compute the list $S= (s_1, \dotsc, s_r)$ of the signatures of $\varphi$ under all real embeddings $\rho_1, \dotsc, \rho_r$ of $K$ and take the maximum of the absolute values of these signatures 
\[
N:=\max_{1\leq j\leq r} |\sgn s_j|.
\]
\item If $N\geq 3$, then return $N$ and quit.
\item Let $\CL$ be the set consisting of all even primes of $K$ and all odd primes dividing any of the coefficients of $\varphi$. 
\item For every $\gp\in \CL$ compute the dimension $d_\gp$ of the anisotropic part of $\varphi\otimes K_\gp$ using Algorithm~\ref{alg_aniso_part_local_dim} and let $M = \max\bigl\{ d_\gp\st \gp\in \CL\bigr\}$.
\item Return $\max\{ M, N\}$.
\end{enumerate}
\end{alg}

\begin{poc}
Let $\psi$ be the anisotropic part of $\varphi$. Obviously 
\[ 
\dim\psi\equiv \dim\varphi\equiv |\sgn \rho_j(\varphi)|\pmod{2}
\] 
for any real embedding $\rho_j$ of $K$. Take $N$ to be the maximum of the absolute values of the signatures of $\varphi$ at all the real places. Now, $\psi$ being anisotropic must be anisotropic at some place of $K$, either finite or infinite. Therefore, clearly $\dim \psi$ is the maximum of the dimensions of the anisotropic parts of the localizations of $\varphi$ at all the places of $K$. However, if $N\geq 3$, then we do not need to consider the finite primes at all. Indeed, if $\dim\psi\geq 5$, then \cite[Theorem~VI.2.2]{Lam05} implies that it must be an infinite, hence real, place. Therefore in this case $\dim\psi= N\geq 5$. Similarly, if $N=3$ or $N=4$, then there is a real embedding $\rho_j$, with $\sgn \rho_j(\psi)= \sgn \rho_j(\varphi)= N$ and so $\dim\psi\geq N$. On the other hand, it cannot be strictly greater, since otherwise $\psi$ would have to be anisotropic at some finite place contrary to the already mentioned  \cite[Theorem~VI.2.2]{Lam05}.
\end{poc}

\section{Level of a number field}\label{sec_level}
In this section we present an algorithm determining an important invariant of a number field, namely its level. Recall that a \term{level} of a field $K$, denoted $s(K)$ is the minimal number of terms needed to represent $-1$ as a sum of squares in $K$. We set $s(K)= \infty$, when $-1$ cannot be expressed as a sum of squares (i.e. $K$ is formally real). 

\begin{alg}\label{alg_level}
Given a number field $K=\QQt$ specified by its defining polynomial $f$, this algorithm computes the level $s(K)$.
\begin{enumerate}
\item If the degree of $f$ is odd, then return $\infty$ and quit.
\item\label{st_is_real} Check if $f$ has any real roots (see Remark~\ref{rem_root_counting} below). If so, then return $\infty$ and quit.
\item Check if $-1$ is a square in $K$. If so, then return $1$ and quit.
\item Find the factorization of $2$ in $\CO_K$ in the form of a list $\CL$ consisting of triples $(\gd_j, e_j, f_j)$, where $\gd_j$ is a prime of $K$ dominating $2$ with the ramification index $e_j$ and the inertia degree $f_j$.
\item\label{step_odd_eifi} If for any $j$, both $e_j$ and $f_j$ are odd, then return $4$ and quit.
\item Return $2$.
\end{enumerate}
\end{alg}

\begin{poc}
The real roots of $f$ correspond to real embeddings of $K$. Hence, if $f$ has a real root (this happens trivially, when $\deg f$ is odd), then $K$ is formally real and consequently its level equals $s(K) = \infty$. Next, if $-1$ is a square in $K$, then $s(K)= 1$. In every other case, $s(K)$ is either $2$ or $4$. In order to distinguish between these two cases, \cite[Proposition XI.2.11]{Lam05} comes in handy. It asserts that $s(K) = 4$ if and only if there is $1\leq j\leq k$ such that $d_j = e_jf_j$ is odd. Otherwise, $s(K) = 2$. This is precisely what step~\eqref{step_odd_eifi} at the end of the algorithm is~for.
\end{poc}

\section{Pythagoras number}
Another field invariant, important from the point of view of the algebraic theory of quadratic forms, is the Pythagoras number. It turns out that an algorithm computing it is in principle the same as the one computing the level. Recall (see e.g. \cite[Chapter~XI]{Lam05}) that a Pythagoras number $P(K)$ of a field $K$ is the smallest integer $p\in \NN$ such that every sum of squares in $K$ is a sum of $p$ squares. If no such an integer exists, then $P(K):= \infty$. It is well known (see e.g. \cite[Theorem~XI.5.6]{Lam05}) that for any arbitrary non-real field $K$ (not necessarily a number field), its Pythagoras number $P(K)$ and its level $s(K)$ differ by no more than $1$, more precisely
\begin{equation}\label{eq_level_and_Pythagoras}
P(K) = s(K)\quad\text{or}\quad P(K) = s(K)+1.
\end{equation}
There are known examples of fields for which any of these two equalities holds. Nevertheless, for number fields the situation is much simpler. We claim that for number fields the latter case is possible only when $s(K) = 4$. We expect that the following result is known to the experts in the field, but since we are not aware of any easily available reference, it easier just to prove it.

\begin{prop}\label{prop_s(K)->P(K)}
Let $K$ be a non-real number field, then
\[
P(K) = 
\begin{cases}
2, & \text{if }s(K) = 1\\
3, & \text{if }s(K) = 2\\
4, & \text{if }s(K) = 4.
\end{cases}
\]
\end{prop}

\begin{proof}
The quadratic closure of $\QQ$ has infinite degree over $\QQ$, hence it cannot be contained in any algebraic number field. It follows that $P(K)\neq 1$ for any number field $K$. On the other hand, $P(K)\leq s(K)+1$ by Eq.~\eqref{eq_level_and_Pythagoras}. Therefore $P(K)=2$ whenever $s(K) = 1$.

Now, assume that $s(K) = 2$, we need to show that $P(K) = 3$. Suppose otherwise, i.e. suppose that $P(K) = s(K)= 2$. This means that for every prime $\gp$ of $K$ and every element $a\in \dot{K}$, the form $\form{a,1,1}\otimes K_\gp$ is isotropic. Fix first an odd prime $\gp$ and take $a$ to be its uniformizer. Then $\form{a,1,1}\otimes K_\gp$ is isotropic if and only if $\form{1,1}\otimes K_\gp$ is isotropic (by \cite[Proposition~VI.1.9]{Lam05}) and so $-1$ is a square in~$K_\gp$. Now take an even prime $\gd$. By our assumption, for every $a\in \dot{K}$, the form $\varphi_\gd= \form{a,1,1}\otimes K_\gd$ is isotropic. It follows form \cite[Proposition~V.3.22]{Lam05} that the Hasse invariant $s(\varphi_\gd)$ of $\varphi_\gd$ equals $(-1,-\det\varphi_\gd)_\gd = (-1,-1)_\gd$. Now the Hasse invariant of $\varphi_\gd$ is
\[
s(\varphi_\gd) = (a,1)_\gd(a,1)_\gd(1,1)_\gd= 1.
\]
Thus, $(-1,-1)_\gd= 1$ for every $a\in \dot{K}$. By the non-degeneracy  of the Hilbert symbol (see e.g. \cite[Theorem~VI.2.16]{Lam05}), this means that $-1$ is a square in $K_\gd$. All in all, we showed that $-1$ is a square in every completion of $K$, but then $s(K)= 1$ contrary to our assumption $s(K) = 2$. This proves the claim $P(K) = 3$.

Finally assume that $s(K) = 4$. A form $\form{a,1,1,1,1}$ is isotropic over $K$ for every $a\in \dot{K}$ by \cite[Corollary~VI.3.5]{Lam05}. Hence, every $a\in \dot{K}$ is a sum of four squares and consequently $P(K)= 4$.
\end{proof}

\begin{alg}\label{alg_Pythagoras_number}
Given a number field $K=\QQt$ specified by its defining polynomial $f$, this algorithm computes the Pythagoras number $P(K)$.
\begin{enumerate}
\item Check if $-1$ is a square in $K$. If so, then return $P(K) = 2$ and quit.
\item Find the factorization of $2$ in $\CO_K$ in the form of a list $\CL$ consisting of triples $(\gd_j, e_j, f_j)$, where $\gd_j$ is a prime of $K$ dominating $2$ with the ramification index $e_j$ and the inertia degree $f_j$.
\item\label{step_P(K)_odd_eifi} If for any $j$, both $e_j$ and $f_j$ are odd, then return $P(K) = 4$ and quit.
\item Return $3$.
\end{enumerate}
\end{alg}

\begin{poc}
Fix a number field $K$. If $K$ is formally real, then \cite[Example~XI.5.9]{Lam05} asserts that $P(K) = 4$ iff there is an even prime $\gd_j$ such that $(K_{\gd_j}: \QQ_2)$ is odd, otherwise $P(K) = 3$. Now, $(K_{\gd_j}: \QQ_2) = e_j f_j$ and so the above-mentioned condition is equivalent to the test in step~\eqref{step_P(K)_odd_eifi}. On the other hand, if $K$ is not real, then the correctness of the algorithm follows immediately from \cite[Proposition~XI.2.11]{Lam05} and Proposition~\ref{prop_s(K)->P(K)}.
\end{poc}

\section{Witt equivalence}\label{sec_Witt_eq}
An important problem in the algebraic theory of quadratic forms is to find criteria for an existence of an isomorphism between the Witt rings of two fields. Such fields are then called \term{Witt equivalent} if the above-mentioned isomorphism exists. In this section we present an algorithm computing the complete set of Witt equivalence invariants of a given number fields. In particular, comparing the results returned by the algorithm one can check whether two number fields are Witt equivalent or not. It was proved in~\cite{Szymiczek91} that the following invariants fully determine the Witt class of a number field $K$:
\begin{itemize}
\item $d = ( K\colon \QQ )$ the degree of $K$ over $\QQ$;
\item $r$ the number of real embeddings of $K$;
\item $s = s(K)$ the level of $K$;
\item $k$ the number of dyadic primes of $K$;
\item for each dyadic prime $\gd_j$ with $1\leq j\leq k$, the pair $(d_j, s_j)$ consisting of a local degree $d_j = (K_{\gd_j}\colon \QQ_2)$ and the local level $s_j = s(K_{\gd_j})$.
\end{itemize}
We claim that all these invariants are computable. 

Let again $K=\QQt$ be a fixed number field specified by the minimal polynomial $f\in \QQ[x]$ of the generator~$\vartheta$. The first two invariants $d$ and $r$ are trivially computable. The degree $d$ is just the degree $\deg f$ of the defining polynomial. In order to compute $r$ one simply counts the number of real roots of $f$ (see Remark~\ref{rem_root_counting} below). In the previous section we showed how to compute the level of $K$. This leaves us only with the local invariants. Assume that the principal ideal $2\CO_K$ factors into prime ideals as:
\[
2\CO_K = \gd_1^{e_1}\dotsb \gd_k^{e_k}
\]
and let $f_j= (\CO_K/\gd_j \colon \FF_2)$ be the inertia degree of $\gd_j$ ($1\leq j\leq k)$. The local degree $d_j = (K_{\gd_j}\colon \QQ_2)$ is the product $d_j = e_j f_j$. What we need is to determine the local level $s_j = s(K_{\gd_j})$. Fix an even prime $\gd = \gd_j$.


\begin{alg}\label{alg_dyadic_level}
Let $\gd$ be an even prime of a number field $K$, $e$ be the ramification index and $f$ the inertia degree of $\gd$. This algorithm computes the level $s(K_\gd)$ of the dyadic completion $K_\gd$ of~$K$.
\begin{enumerate}
\item Check if $-1$ is a square in $K$, if so then return $1$ and quit.
\item If both $e$ and $f$ are odd, then return $4$.
\item If $e$ is odd but $f$ is even, then return $2$.
\item If $e$ is even, check whether $-1$ is a square in $K_\gd$, if so then return $1$, if not then return $2$.
\end{enumerate}
\end{alg}

\begin{poc}
It is clear that if $-1$ is a square already in $K$, then it is also a square in $K_\gd$ and so $s(K_\gd)=1$. This justifies the first step. Suppose that $e$ is odd. Let $\QQ_2(\eta)$ be the (unique) maximal unramified extension of $\QQ_2$ contained in~$K_{\gd}$. Since the quadratic extension $\QQ_2(i)/\QQ_2$ is totally ramified (see e.g. \cite[Ch.~V \S2]{Narkiewicz90}), it follows that $i\notin \QQ_2(\eta)$. Now, $\bigl(\QQ_2(\eta)\colon \QQ_2\bigr) = f$ and $(K_{\gd}\colon \QQ_2)= e f$. Hence the relative degree $\bigl(K_{\gd}\colon \QQ_2(\eta)\bigr)$ equals $e$ and so is odd. In particular $i\notin K_{\gd}$ and so $s(K_\gd)\geq 2$. Finally \cite[Example~XI.2.4]{Lam05} asserts that $s(K_\gd)=4$ if and only if $(K_{\gd}\colon \QQ_2)$ is odd.

Conversely, assume that $e$ is even and so is the degree $(K_{\gd}\colon \QQ_2)$. It follows from \cite[Example~XI.2.4]{Lam05} that $s(K_\gd)\leq 2$. It equals one if and only if $-1$ is a square in~$K_\gd$.
\end{poc}

Having all the necessary ingredients ready we may now present the last algorithm of this paper that construct the complete set of Witt equivalence invariants.

\begin{alg}\label{alg_WE_invariants}
If $K = \QQt$ is a number field specified \textup(up to an isomorphism\textup) by the minimal polynomial $f\in \QQ[t]$ of its generator, then this algorithm computes the complete set of Witt equivalence invariants of $K$. In particular, two fields are Witt equivalent if and only if the outputs of the algorithm are the same for both fields.
\begin{enumerate}
\item Let $d = \deg f$.
\item\label{st_count_r_roots} Compute the number $r$ of real roots of $f$.
\item Use Algorithm~\ref{alg_level} to compute the level $s= s(K)$.
\item Let $\CL = \bigl\{ (\gd_j, e_j, f_j) \bigr\}$ be the list of all even primes of $K$ together with their ramification indices and inertia degrees. 
\item Take an empty list $\CS$.
\item For each even prime $\gd_j\in \CL$ let $d_j = e_j f_j$. Use Algorithm~\ref{alg_dyadic_level} to compute the local level $s_j = s(K_{\gd_j})$. Append the pair $(d_j, s_j)$ to the list $\CS$.
\item Sort the list $\CS$ lexicographically.
\item Return $(d,r,s,k,\CS)$.
\end{enumerate}
\end{alg}

\begin{rem}\label{rem_root_counting}
There is a number of known algorithms which can be used to count real roots of~$f$ in step~\eqref{st_count_r_roots}. They vary from methods based on Sturm's and Hermite's theorems (see e.g. \cite[Theorems~2.56, 4.13 and also Algorithm~9.28]{BPR03}) to those based on Vincent's theorem (see \cite{AS05, AV10}). Of course, any algorithm that counts real roots can also be used to check if a polynomial has at least one real root, which is needed in step~\eqref{st_is_real} of Algorithm~\ref{alg_level}. Honestly, the authors of this paper are not aware of any method answering the latter question, that would be significantly simpler than a general root counting algorithm.
\end{rem}

\section{Example applications}\label{sec_examples}
In order to verify the correctness of the algorithm as well as to allow experimentation, we implemented the presented algorithm in a computer algebra system Sage (see \cite{sage}). The code is available from the first author's home page at \url{http://z2.math.us.edu.pl/perry/papersen.html}. A formula for the number of Witt classes of number fields of a fixed degree was developed in \cite{Szymiczek91}. Nevertheless, actual representatives for these classes were only found for quadratic and cubic fields in \cite{Szymiczek91} and for quartic fields in \cite{JakMar92}. The first test for usability of our implementation was to find new representatives of all classes of cubic and quartic fields. Next, we found the representatives of all $36$ classes of quintic fields  and all $95$ classes of sextic fields. These two results are completely new. The findings are gathered in tables \ref{tbl_cubic}--\ref{tbl_sextic}. For those classes, for which we found more than one field, the corresponding table contains a representative with the smallest absolute value of the discriminant.

The method used here was a combination of an `aided random search' (explained bellow) and (following a suggestion of the reviewer) a search of the data base of number fields in~\cite{lmfdb}. In case of quartic fields we were able to significantly improve the known results, as the largest discriminant in our case is $122\,825$ and the largest absolute value of a coefficient of a defining polynomial is $86$ (vs. respectively $210\,668\,284$ and $208\,042$ in \cite{JakMar92}).

Some of Witt classes are extremely rare and virtually impossible to be found by a blind random search. These are mostly the classes of fields were $2$ splits completely. In order to find the representatives of these classes one may proceed as follows. Denote by $|a|_2:= 2^{-\ord_2 a}$ the canonical dyadic norm and let $\| (a_0, \dotsc, a_d)\|_2:= \max \bigl\{ |a_0|_2, \dotsc, |a_d|_d\bigr\}$ be the associated norm of the vector space $(\QQ_2)^{d+1}$. Take a polynomial $f$ with $d$ distinct integral roots. Write it as a dot product $f = V\cdot X$, where $V$ is the vector of coefficients of $f$ and $X= (1,x, x^2, \dotsc, x^d)^T$ are the powers of $x$. Take now some random vector $W$ and let $\tilde{f} = (V+W)\cdot X$. If the norm $\|W\|_2$ of $W$ is small enough, then $\tilde{f}$ still has $d$ distinct roots in $\QQ_2$ but there is a good chance that it is irreducible over $\QQ$. It follows that $2$ splits completely in the field $K = \QQ[x]/\form{ \tilde{f}}$, as desired.


\subsection{Conner's Problem} 
A number field $K$ is said to satisfy Conner's Level Condition (CLC for short) if $s(K) = 2$ but $s(K_\gd) = 1$ for every even prime $\gd$ of $K$. \cite{JMSz95, JMSz97} proved that if a number field satisfies CLC, then its class number is even.

Since CLC is expressed in terms of Witt equivalence invariants, thus one may treat it as a property of Witt equivalence classes. In particular, if a Witt equivalence class satisfies CLC, then every field in this class has an even class number. P.E.~Conner asked for an inverse of this statement (c.f. \cite{Szymiczek2000}):\smallskip
\begin{center}
\emph{Suppose a Witt equivalence class does not satisfy CLC.\\ Does it contain a field with an odd class number?}
\end{center}\smallskip
An affirmative answer to Conner's question was found in \cite{Szymiczek91} for quadratic and cubic fields and in \cite{JMSz95, JMSz97} for quartic fields.

Observe that a field of an odd degree cannot satisfy CLC (since its level is infinite). Using \cite{lmfdb} one checks that all fields in Tables~\ref{tbl_cubic} and~\ref{tbl_quintic} have trivial class groups.

As for quartic and sextic fields, there are precisely five Witt equivalence classes for which CLC holds. These are classes: \ref{c4_CLC_1}, \ref{c4_CLC_2}, \ref{c6_CLC_1}, \ref{c6_CLC_2} and~\ref{c6_CLC_3}. Once we omit them, all other representatives listed in tables~\ref{tbl_quartic} and~\ref{tbl_sextic}, except \ref{c4_h=2} (which has class number~$2$), have odd class numbers. In fact, all these fields have class numbers not only odd but actually equal one, with only two exceptions. The exceptions are the representatives of Witt equivalence classes: \ref{c6_h=9} and~\ref{c6_h=5}. Their class numbers equal respectively: $9$ and~$5$. Nevertheless, it is possible to find representatives of these three ``exceptional'' classes with trivial class groups but with higher absolute values of discriminants (recall that the tables contain representatives with smallest discriminant we were able to find):\smallskip
\begin{center}
\begin{tabular}{|r|l|l|}
\hline
class & defining polynomial & LMFDB label \\\hline
\ref{c4_h=2} & $x^4 - 2x^3 - 13x^2 + 14x + 32$ & 4.4.164441.1\\
\ref{c6_h=9} & $x^6 - 3x^5 - 21x^4 - x^3 + 228x^2 + 532x + 448$ & 6.0.827250487.1\\
\ref{c6_h=5} & $x^6 + 2x^4 + x^2 + 28$ & 6.0.12122992.1\\\hline
\end{tabular}
\end{center}\medskip
In all cases, the class numbers where either obtained from \cite{lmfdb} or computed in Sage using GP/Pari back-end (see \cite{pari}), except for classes \ref{c6_GRH_28}, \ref{c6_GRH_50}--\ref{c6_GRH_52}, \ref{c6_GRH_54}, \ref{c6_GRH_74}--\ref{c6_GRH_78}. The class numbers of these ten fields were computed using Magma back-end (see \cite{magma}) under assumption of Generalized Riemann Hypothesis.

Summarizing the above discussion, this proves:

\begin{thm}
The Conner's question has an affirmative answer for fields of degree $<7$. What is more, every Witt equivalence class of number fields of degree $<7$ that do not satisfy CLC contains a field with a trivial class group.
\end{thm}



\FloatBarrier\appendix
\section{Tables of representatives of Witt classes}
\newcounter{LP}
\newcounter{deg}
\renewcommand{\theLP}{\arabic{deg}.\arabic{LP}}
\newcommand{\lp}{\refstepcounter{LP}\theLP.}
\renewcommand\arraystretch{1.3}

\FloatBarrier\setcounter{LP}{0}\setcounter{deg}{3}
\begin{longtable}{|r|>{\begin{flushleft}}p{36mm}<{\end{flushleft}}|l|c|c|}
\caption{\label{tbl_cubic}Witt classes of cubic fields}\\
\hline
No. & defining polynomial & LMFDB & $r$ & dyadic degrees and levels\\
\hline\endfirsthead
\caption{Witt classes of cubic fields}\\
\hline
No. & defining polynomial & LMFDB & $r$ & dyadic degrees and levels\\
\hline\endhead
\hline \multicolumn{5}{r}{\textit{Continued on the next page}} \\
\endfoot
\hline
\endlastfoot
\lp & $ x^{3} - x - 8 $           &  3.1.431.1 & $1$ & $\bigl\{ (1, 4), (1, 4), (1, 4) \bigr\}$\\
\lp & $ x^{3} + 2x - 1 $          &  3.1.59.1  & $1$ & $\bigl\{ (1, 4), (2, 1) \bigr\}$\\
\lp & $ x^{3} - 3x - 4 $          &  3.1.324.1 & $1$ & $\bigl\{ (1, 4), (2, 2) \bigr\}$\\
\lp & $ x^{3} - x^{2} + 1 $       &  3.1.23.1  & $1$ & $\bigl\{ (3, 4) \bigr\}$\\
\lp & $ x^{3} - x^{2} - 10x + 8 $ &  3.3.961.1 & $3$ & $\bigl\{ (1, 4), (1, 4), (1, 4) \bigr\}$\\
\lp & $ x^{3} - 4x - 1 $          &  3.3.229.1 & $3$ & $\bigl\{ (1, 4), (2, 1) \bigr\}$\\
\lp & $ x^{3} - x^{2} - 4x + 2 $  & 3.3.316.1  & $3$ & $\bigl\{ (1, 4), (2, 2) \bigr\}$\\
\lp & $ x^{3} - x^{2} - 2x + 1 $  & 3.3.49.1   & $3$ & $\bigl\{ (3, 4) \bigr\}$\\
\end{longtable}

\FloatBarrier\setcounter{LP}{0}\setcounter{deg}{4}
\newcommand{\SSS}[1]{\bigl\{#1\bigr\}}
\begin{longtable}{|r|>{\begin{flushleft}}p{36mm}<{\end{flushleft}}|l|c|c|c|}
\caption{\label{tbl_quartic}Witt classes of quartic fields}\\
\hline
No. & defining polynomial & LMFDB & $r$ & $s$ & dyadic degrees and levels\\
\hline\endfirsthead
\caption{Witt classes of quartic fields (continued)}\\
\hline
No. & defining polynomial & LMFDB  & $r$ & $s$ & dyadic degrees and levels\\
\hline\endhead
\hline \multicolumn{6}{r}{\textit{Continued on the next page}} \\
\endfoot
\hline
\endlastfoot
\lp & $ x^{4} - 2x^{3} - x^{2} + 2x + 2 $ & 4.0.656.1 & $ 0 $ & $ 1 $ & $\bigl\{ (2, 1), (2, 1) \bigr\}$\\
\lp & $ x^{4} - x^{2} + 1 $ & 4.0.144.1 & $ 0 $ & $ 1 $ & $\bigl\{ (4, 1) \bigr\}$\\
\lp\label{c4_CLC_1} & $ x^{4} + 3x^{2} - 14x + 18 $ & 4.0.44688.1 & $ 0 $ & $ 2 $ & $\bigl\{ (2, 1), (2, 1) \bigr\}$\\
\lp & $ x^{4} - x^{3} + x^{2} + 4x + 2 $ & 4.0.2156.1 & $ 0 $ & $ 2 $ & $\bigl\{ (2, 1), (2, 2) \bigr\}$\\
\lp & $ x^{4} - x^{3} + 2x^{2} + x + 1 $ & 4.0.225.1 & $ 0 $ & $ 2 $ & $\bigl\{ (2, 2), (2, 2) \bigr\}$\\
\lp\label{c4_CLC_2} & $ x^{4} - 5x^{2} + 25 $ & 4.0.3600.3 & $ 0 $ & $ 2 $ & $\bigl\{ (4, 1) \bigr\}$\\
\lp & $ x^{4} - x^{3} - x^{2} + x + 1 $ & 4.0.117.1 & $ 0 $ & $ 2 $ & $\bigl\{ (4, 2) \bigr\}$\\
\lp & $ x^{4} - 2x^{3} - x^{2} + 2x + 8 $ & 4.0.6713.1 & $ 0 $ & $ 4 $ & $\bigl\{ (1, 4), (1, 4), (1, 4), (1, 4) \bigr\}$\\
\lp & $ x^{4} - x^{3} + 6x^{2} - 2x + 4 $ & 4.0.4508.1 & $ 0 $ & $ 4 $ & $\bigl\{ (1, 4), (1, 4), (2, 1) \bigr\}$\\
\lp & $ x^{4} - 2x^{3} + 2x^{2} - x + 2 $ & 4.0.1421.1 & $ 0 $ & $ 4 $ & $\bigl\{ (1, 4), (1, 4), (2, 2) \bigr\}$\\
\lp & $ x^{4} + x^{2} - x + 1 $ & 4.0.257.1 & $ 0 $ & $ 4 $ & $\bigl\{ (1, 4), (3, 4) \bigr\}$\\
\lp & $ x^{4} - 2x^{3} - 5x^{2} + 6x - 8 $ & 4.2.29767.1 & $ 2 $ & $ \infty $ & $\bigl\{ (1, 4), (1, 4), (1, 4), (1, 4) \bigr\}$\\
\lp & $ x^{4} - 5x^{2} - 4 $ & 4.2.6724.1 & $ 2 $ & $ \infty $ & $\bigl\{ (1, 4), (1, 4), (2, 1) \bigr\}$\\
\lp & $ x^{4} - 2x^{2} - x - 2 $ & 4.2.4027.1 & $ 2 $ & $ \infty $ & $\bigl\{ (1, 4), (1, 4), (2, 2) \bigr\}$\\
\lp & $ x^{4} - 2x^{3} + x^{2} - x - 1 $ & 4.2.751.1 & $ 2 $ & $ \infty $ & $\bigl\{ (1, 4), (3, 4) \bigr\}$\\
\lp & $ x^{4} - 2x^{3} - 5x^{2} - 2x + 2 $ & 4.2.27632.1 & $ 2 $ & $ \infty $ & $\bigl\{ (2, 1), (2, 1) \bigr\}$\\
\lp & $ x^{4} - x^{3} - 3x^{2} + 2 $ & 4.2.1588.1 & $ 2 $ & $ \infty $ & $\bigl\{ (2, 1), (2, 2) \bigr\}$\\
\lp & $ x^{4} - x^{3} - 3x - 1 $ & 4.2.775.1 & $ 2 $ & $ \infty $ & $\bigl\{ (2, 2), (2, 2) \bigr\}$\\
\lp & $ x^{4} + x^{2} - 6x + 1 $ & 4.2.3312.2 & $ 2 $ & $ \infty $ & $\bigl\{ (4, 1) \bigr\}$\\
\lp & $ x^{4} - x^{3} + 2x - 1 $ & 4.2.275.1 & $ 2 $ & $ \infty $ & $\bigl\{ (4, 2) \bigr\}$\\
\lp\label{c4_h=2} & $ x^{4} - x^{3} - 23x^{2} + x + 86 $ & 4.4.122825.1 & $ 4 $ & $ \infty $ & $\bigl\{ (1, 4), (1, 4), (1, 4), (1, 4) \bigr\}$\\
\lp & $ x^{4} - x^{3} - 15x^{2} + 31x - 8 $ & 4.4.54332.1 & $ 4 $ & $ \infty $ & $\bigl\{ (1, 4), (1, 4), (2, 1) \bigr\}$\\
\lp & $ x^{4} - 2x^{3} - 4x^{2} + 5x + 2 $ & 4.4.15317.1 & $ 4 $ & $ \infty $ & $\bigl\{ (1, 4), (1, 4), (2, 2) \bigr\}$\\
\lp & $ x^{4} - x^{3} - 4x^{2} + x + 2 $ & 4.4.2777.1 & $ 4 $ & $ \infty $ & $\bigl\{ (1, 4), (3, 4) \bigr\}$\\
\lp & $ x^{4} - 2x^{3} - 5x^{2} + 6x + 2 $ & 4.4.44688.2 & $ 4 $ & $ \infty $ & $\bigl\{ (2, 1), (2, 1) \bigr\}$\\
\lp & $ x^{4} - 5x^{2} + 2 $ & 4.4.9248.1 & $ 4 $ & $ \infty $ & $\bigl\{ (2, 1), (2, 2) \bigr\}$\\
\lp & $ x^{4} - x^{3} - 5x^{2} + 2x + 4 $ & 4.4.2225.1 & $ 4 $ & $ \infty $ & $\bigl\{ (2, 2), (2, 2) \bigr\}$\\
\lp & $ x^{4} - 2x^{3} - 7x^{2} + 8x + 1 $ & 4.4.3600.1 & $ 4 $ & $ \infty $ & $\bigl\{ (4, 1) \bigr\}$\\
\lp & $ x^{4} - x^{3} - 3x^{2} + x + 1 $ & 4.4.725.1 & $ 4 $ & $ \infty $ & $\bigl\{ (4, 2) \bigr\}$\\
\hline
\end{longtable}

\FloatBarrier\setcounter{LP}{0}\setcounter{deg}{5}
\begin{landscape}
\begin{longtable}{|r|>{\begin{flushleft}}p{65mm}<{\end{flushleft}}|l|c|c|}
\caption{\label{tbl_quintic}Witt classes of quintic fields}\\
\hline
No. & defining polynomial & LMFDB label & $r$ & dyadic degrees and levels\\
\hline\endfirsthead
\caption{Witt classes of quintic fields (continued)}\\
\hline
No. & defining polynomial & LMFDB label & $r$ & dyadic degrees and levels\\
\hline\endhead
\hline \multicolumn{5}{r}{\textit{Continued on the next page}} \\
\endfoot
\hline
\endlastfoot
\lp & $ x^{5} - x^{3} - 4 \, x^{2} + 20 \, x + 16 $ & 5.1.1659001.1 & $ 1 $ & $\bigl\{ (1, 4), (1, 4), (1, 4), (1, 4), (1, 4) \bigr\}$\\
\lp & $ x^{5} - 2 \, x^{4} + 8 \, x^{3} - 8 \, x^{2} + 15 \, x + 2 $ & 5.1.224956.1 & $ 1 $ & $\bigl\{ (1, 4), (1, 4), (1, 4), (2, 1) \bigr\}$\\
\lp & $ x^{5} - 3 \, x^{2} + 8 \, x - 4 $ & 5.1.119701.2 & $ 1 $ & $\bigl\{ (1, 4), (1, 4), (1, 4), (2, 2) \bigr\}$\\
\lp & $ x^{5} - 2 \, x^{4} + 3 \, x^{3} - 5 \, x^{2} + 3 \, x - 2 $ & 5.1.19633.1 & $ 1 $ & $\bigl\{ (1, 4), (1, 4), (3, 4) \bigr\}$\\
\lp & $ x^{5} - x^{3} - 4 \, x^{2} + 6 \, x - 4 $ & 5.1.408976.1 & $ 1 $ & $\bigl\{ (1, 4), (2, 1), (2, 1) \bigr\}$\\
\lp & $ x^{5} + 2 \, x^{3} - x^{2} + 2 \, x - 2 $ & 5.1.28684.1 & $ 1 $ & $\bigl\{ (1, 4), (2, 1), (2, 2) \bigr\}$\\
\lp & $ x^{5} - x^{4} + 2 \, x^{3} - x^{2} + x + 2 $ & 5.1.17161.1 & $ 1 $ & $\bigl\{ (1, 4), (2, 2), (2, 2) \bigr\}$\\
\lp & $ x^{5} - 2 \, x^{4} + x^{3} + 4 \, x^{2} - 7 \, x + 4 $ & 5.1.42256.1 & $ 1 $ & $\bigl\{ (1, 4), (4, 1) \bigr\}$\\
\lp & $ x^{5} - 2 \, x^{2} - 2 \, x - 1 $ & 5.1.4261.1 & $ 1 $ & $\bigl\{ (1, 4), (4, 2) \bigr\}$\\
\lp & $ x^{5} - x^{4} + 2 \, x^{3} - 3 \, x^{2} + 2 \, x - 2 $ & 5.1.10492.1 & $ 1 $ & $\bigl\{ (2, 1), (3, 4) \bigr\}$\\
\lp & $ x^{5} - x - 1 $ & 5.1.2869.1 & $ 1 $ & $\bigl\{ (2, 2), (3, 4) \bigr\}$\\
\lp & $ x^{5} - x^{3} - x^{2} + x + 1 $ & 5.1.1609.1 & $ 1 $ & $\bigl\{ (5, 4) \bigr\}$\\
\lp & $ x^{5} - x^{4} - 3 \, x^{3} + 17 \, x^{2} - 38 \, x + 8 $ & 5.3.6556247.1 & $ 3 $ & $\bigl\{ (1, 4), (1, 4), (1, 4), (1, 4), (1, 4) \bigr\}$\\
\lp & $ x^{5} - 8 \, x^{3} + 7 \, x - 16 $ & 5.3.1078244.1 & $ 3 $ & $\bigl\{ (1, 4), (1, 4), (1, 4), (2, 1) \bigr\}$\\
\lp & $ x^{5} - 2 \, x^{4} + 9 \, x^{2} - 4 $ & 5.3.443291.1 & $ 3 $ & $\bigl\{ (1, 4), (1, 4), (1, 4), (2, 2) \bigr\}$\\
\lp & $ x^{5} - 2 \, x^{4} + x^{3} + x^{2} - 7 \, x + 4 $ & 5.373607.1 & $ 3 $ & $\bigl\{ (1, 4), (1, 4), (3, 4) \bigr\}$\\
\lp & $ x^{5} - x^{3} - 4 \, x^{2} - 2 \, x + 4 $ & 5.3.409328.1 & $ 3 $ & $\bigl\{ (1, 4), (2, 1), (2, 1) \bigr\}$\\
\lp & $ x^{5} - x^{4} - 4 \, x^{3} + x^{2} + 3 \, x + 2 $ & 5.3.113684.1 & $ 3 $ & $\bigl\{ (1, 4), (2, 1), (2, 2) \bigr\}$\\
\lp & $ x^{5} - 4 \, x^{3} - x^{2} + 2 \, x + 4 $ & 5.3.67943.1 & $ 3 $ & $\bigl\{ (1, 4), (2, 2), (2, 2) \bigr\}$\\
\lp & $ x^{5} - 2 \, x^{3} - 2 \, x^{2} - 3 \, x + 2 $ & 5.3.41456.1 & $ 3 $ & $\bigl\{ (1, 4), (4, 1) \bigr\}$\\
\lp & $ x^{5} - x^{4} - x^{3} - x^{2} - 3 \, x + 1 $ & 5.3.13523.1 & $ 3 $ & $\bigl\{ (1, 4), (4, 2) \bigr\}$\\
\lp & $ x^{5} - x^{4} - 2 \, x^{3} - x^{2} + 2 \, x + 2 $ & 5.3.17348.1 & $ 3 $ & $\bigl\{ (2, 1), (3, 4) \bigr\}$\\
\lp & $ x^{5} - x^{4} - 2 \, x + 1 $ & 5.3.11243.1 & $ 3 $ & $\bigl\{ (2, 2), (3, 4) \bigr\}$\\
\lp & $ x^{5} - x^{3} - 2 \, x^{2} + 1 $ & 5.3.4511.1 & $ 3 $ & $\bigl\{ (5, 4) \bigr\}$\\
\lp & $ x^{5} - x^{4} - 19 \, x^{3} + 17 \, x^{2} + 58 \, x - 40 $ & 5.5.46919377.1 & $ 5 $ & $\bigl\{ (1, 4), (1, 4), (1, 4), (1, 4), (1, 4) \bigr\}$\\
\lp & $ x^{5} - x^{4} - 14 \, x^{3} - 2 \, x^{2} + 28 \, x + 16 $ & 5.5.8048764.1 & $ 5 $ & $\bigl\{ (1, 4), (1, 4), (1, 4), (2, 1) \bigr\}$\\
\lp & $ x^{5} - 13 \, x^{3} - 4 \, x^{2} + 24 \, x + 8 $ & 5.5.3609877.1 & $ 5 $ & $\bigl\{ (1, 4), (1, 4), (1, 4), (2, 2) \bigr\}$\\
\lp & $ x^{5} - 7 \, x^{3} - x^{2} + 11 \, x + 4 $ & 5.5.372289.1 & $ 5 $ & $\bigl\{ (1, 4), (1, 4), (3, 4) \bigr\}$\\
\lp & $ x^{5} - 2 \, x^{4} - 10 \, x^{3} + 14 \, x^{2} + 21 \, x - 16 $ & 5.5.46919377.1 & $ 5 $ & $\bigl\{ (1, 4), (2, 1), (2, 1) \bigr\}$\\
\lp & $ x^{5} - x^{4} - 7 \, x^{3} + 6 \, x^{2} + 6 \, x - 4 $ & 5.5.600268.1 & $ 5 $ & $\bigl\{ (1, 4), (2, 1), (2, 2) \bigr\}$\\
\lp & $ x^{5} - x^{4} - 6 \, x^{3} + 5 \, x^{2} + 7 \, x - 4 $ & 5.5.406264.1 & $ 5 $ & $\bigl\{ (1, 4), (2, 2), (2, 2) \bigr\}$\\
\lp & $ x^{5} - 8 \, x^{3} - 2 \, x^{2} + 5 \, x + 2 $ & 5.5.380224.1 & $ 5 $ & $\bigl\{ (1, 4), (4, 1) \bigr\}$\\
\lp & $ x^{5} - x^{4} - 5 \, x^{3} + 3 \, x^{2} + 5 \, x - 2 $ & 5.5.81509.1 & $ 5 $ & $\bigl\{ (1, 4), (4, 2) \bigr\}$\\
\lp & $ x^{5} - 2 \, x^{4} - 4 \, x^{3} + 5 \, x^{2} + 3 \, x - 1 $ & 5.5.218524.1 & $ 5 $ & $\bigl\{ (2, 1), (3, 4) \bigr\}$\\
\lp & $ x^{5} - x^{4} - 5 \, x^{3} + 4 \, x^{2} + 4 \, x - 1 $ & 5.5.117688.1 & $ 5 $ & $\bigl\{ (2, 2), (3, 4) \bigr\}$\\
\lp & $ x^{5} - x^{4} - 4 \, x^{3} + 3 \, x^{2} + 3 \, x - 1 $ & 5.1.14641.1 & $ 5 $ & $\bigl\{ (5, 4) \bigr\}$\\
\hline
\end{longtable}
\end{landscape}

\FloatBarrier\setcounter{LP}{0}\setcounter{deg}{6}
\begin{landscape}
\begin{longtable}{|r|>{\begin{flushleft}}p{65mm}<{\end{flushleft}}|l|c|c|c|}
\caption{\label{tbl_sextic}Witt classes of sextic fields}\\
\hline
No. & defining polynomial & LMFDB label & $r$ & $s$ & dyadic degrees and levels\\
\hline\endfirsthead
\caption{Witt classes of sextic fields (continued)}\\
\hline
No. & defining polynomial & LMFDB label & $r$ & $s$ & dyadic degrees and levels\\
\hline\endhead
\hline \multicolumn{6}{r}{\textit{Continued on the next page}} \\
\endfoot
\hline
\endlastfoot
\lp & $ x^{6} - 2 \, x^{5} + 3 \, x^{4} - 6 \, x^{3} + 6 \, x^{2} - 8 \, x + 8 $ & 6.0.399424.1 & $ 0 $ & $ 1 $ & $\bigl\{ (2, 1), (2, 1), (2, 1) \bigr\}$\\
\lp & $ x^{6} - 2 \, x^{3} + x^{2} + 2 \, x + 2 $ & 6.0.53824.1 & $ 0 $ & $ 1 $ & $\bigl\{ (2, 1), (4, 1) \bigr\}$\\
\lp & $ x^{6} - x^{4} - 2 \, x^{3} + 2 \, x + 1 $ & 6.0.10816.1 & $ 0 $ & $ 1 $ & $\bigl\{ (6, 1) \bigr\}$\\
\lp\label{c6_CLC_1} & $ x^{6} - 2 \, x^{5} - x^{4} + 10 \, x^{3} + 82 \, x^{2} - 8 \, x + 16 $ & 6.0.45118016.1 & $ 0 $ & $ 2 $ & $\bigl\{ (2, 1), (2, 1), (2, 1) \bigr\}$\\
\lp & $ x^{6} - x^{5} + 3 \, x^{4} - 2 \, x^{3} + 13 \, x^{2} - 17 \, x + 11 $ & 6.0.2787888.1 & $ 0 $ & $ 2 $ & $\bigl\{ (2, 1), (2, 1), (2, 2) \bigr\}$\\
\lp & $ x^{6} - 3 \, x^{5} + 4 \, x^{4} - 3 \, x^{3} + 3 \, x^{2} + 2 $ & 6.0.236708.2 & $ 0 $ & $ 2 $ & $\bigl\{ (2, 1), (2, 2), (2, 2) \bigr\}$\\
\lp\label{c6_CLC_2} & $ x^{6} + 4 \, x^{4} - 8 \, x^{3} + 4 \, x^{2} + 1 $ & 6.0.11999296.1 & $ 0 $ & $ 2 $ & $\bigl\{ (2, 1), (4, 1) \bigr\}$\\
\lp & $ x^{6} + 2 \, x^{4} - x^{3} + 3 \, x^{2} + 2 $ & 6.0.122708.1 & $ 0 $ & $ 2 $ & $\bigl\{ (2, 1), (4, 2) \bigr\}$\\
\lp & $ x^{6} - 2 \, x^{4} - 2 \, x^{3} + 4 \, x^{2} + 2 \, x + 1 $ & 6.0.93987.1 & $ 0 $ & $ 2 $ & $\bigl\{ (2, 2), (2, 2), (2, 2) \bigr\}$\\
\lp & $ x^{6} - 3 \, x^{4} - 2 \, x^{3} + 9 \, x^{2} + 12 \, x + 4 $ & 6.0.314928.2 & $ 0 $ & $ 2 $ & $\bigl\{ (2, 2), (4, 1) \bigr\}$\\
\lp & $ x^{6} - x^{5} - x^{4} + 3 \, x^{3} - 2 \, x + 1 $ & 6.0.16551.1 & $ 0 $ & $ 2 $ & $\bigl\{ (2, 2), (4, 2) \bigr\}$\\
\lp\label{c6_CLC_3} & $ x^{6} - 2 \, x^{5} + 12 \, x^{4} - 30 \, x^{3} + 74 \, x^{2} - 88 \, x + 82 $ & 6.0.9199872.1 & $ 0 $ & $ 2 $ & $\bigl\{ (6, 1) \bigr\}$\\
\lp & $ x^{6} - x^{5} + x^{4} - 2 \, x^{3} + 4 \, x^{2} - 3 \, x + 1 $ & 6.0.9747.1 & $ 0 $ & $ 2 $ & $\bigl\{ (6, 2) \bigr\}$\\
\lp\label{c6_h=9} & $ x^{6} - x^{5} + 3 \, x^{4} - 11 \, x^{3} + 44 \, x^{2} - 36 \, x + 32 $ & 6.0.28629151.1 & $ 0 $ & $ 4 $ & $\bigl\{ (1, 4), (1, 4), (1, 4), (1, 4), (1, 4), (1, 4) \bigr\}$\\
\lp & $ x^{6} - x^{2} + 16 $ & 6.0.11930116.1 & $ 0 $ & $ 4 $ & $\bigl\{ (1, 4), (1, 4), (1, 4), (1, 4), (2, 1) \bigr\}$\\
\lp & $ x^{6} - 3 \, x^{5} + 2 \, x^{4} + x^{3} + x^{2} - 2 \, x + 8 $ & 6.0.4807171.1 & $ 0 $ & $ 4 $ & $\bigl\{ (1, 4), (1, 4), (1, 4), (1, 4), (2, 2) \bigr\}$\\
\lp & $ x^{6} - x^{5} + 3 \, x^{4} + 2 \, x^{3} + 6 \, x^{2} + 3 \, x + 2 $ & 6.0.469567.1 & $ 0 $ & $ 4 $ & $\bigl\{ (1, 4), (1, 4), (1, 4), (3, 4) \bigr\}$\\
\lp\label{c6_h=5} & $ x^{6} - x^{5} + 6 \, x^{4} + 4 \, x^{3} + 11 \, x^{2} + 21 \, x + 22 $ & 6.0.7888624.1 & $ 0 $ & $ 4 $ & $\bigl\{ (1, 4), (1, 4), (2, 1), (2, 1) \bigr\}$\\
\lp & $ x^{6} - 13 \, x^{2} + 20 $ & 6.0.5060180.2 & $ 0 $ & $ 4 $ & $\bigl\{ (1, 4), (1, 4), (2, 1), (2, 2) \bigr\}$\\
\lp & $ x^{6} - 3 \, x^{5} + x^{4} + 3 \, x^{3} + x^{2} - 3 \, x + 2 $ & 6.0.489119.1 & $ 0 $ & $ 4 $ & $\bigl\{ (1, 4), (1, 4), (2, 2), (2, 2) \bigr\}$\\
\lp & $ x^{6} - x^{5} - x^{4} - 3 \, x^{3} + x^{2} + 7 \, x + 4 $ & 6.0.887152.1 & $ 0 $ & $ 4 $ & $\bigl\{ (1, 4), (1, 4), (4, 1) \bigr\}$\\
\lp & $ x^{6} - 2 \, x^{5} + 2 \, x^{4} - x + 2 $ & 6.0.134363.1 & $ 0 $ & $ 4 $ & $\bigl\{ (1, 4), (1, 4), (4, 2) \bigr\}$\\
\lp & $ x^{6} - 2 \, x^{5} + x^{4} - x^{3} + x^{2} + 2 $ & 6.0.1085252.1 & $ 0 $ & $ 4 $ & $\bigl\{ (1, 4), (2, 1), (3, 4) \bigr\}$\\
\lp & $ x^{6} + x^{4} - x + 1 $ & 6.0.89363.1 & $ 0 $ & $ 4 $ & $\bigl\{ (1, 4), (2, 2), (3, 4) \bigr\}$\\
\lp & $ x^{6} - x^{5} + x^{4} - x^{3} + 2 \, x^{2} - x + 1 $ & 6.0.31223.1 & $ 0 $ & $ 4 $ & $\bigl\{ (1, 4), (5, 4) \bigr\}$\\
\lp & $ x^{6} - 3 \, x^{5} + 5 \, x^{4} - 5 \, x^{3} + 5 \, x^{2} - 3 \, x + 1 $ & 6.0.12167.1 & $ 0 $ & $ 4 $ & $\bigl\{ (3, 4), (3, 4) \bigr\}$\\
\lp & $ x^{6} - 53 \, x^{4} - 16 \, x^{3} + 868 \, x^{2} - 800 \, x - 4288 $ & --- & $ 2 $ & $ \infty $ & $\bigl\{ (1, 4), (1, 4), (1, 4), (1, 4), (1, 4), (1, 4) \bigr\}$\\
\lp\label{c6_GRH_28} & $ x^{6} - 232 \, x^{5} - 479 \, x^{4} - 440 \, x^{3} - 502 \, x^{2} + 348 \, x - 64 $ & --- & $ 2 $ & $ \infty $ & $\bigl\{ (1, 4), (1, 4), (1, 4), (1, 4), (2, 1) \bigr\}$\\
\lp & $ x^{6} - 3 \, x^{5} - 4 \, x^{4} + 13 \, x^{3} - x^{2} - 6 \, x - 8 $ & 6.2.26919373.1 & $ 2 $ & $ \infty $ & $\bigl\{ (1, 4), (1, 4), (1, 4), (1, 4), (2, 2) \bigr\}$\\
\lp & $ x^{6} - x^{5} - x^{4} - 7 \, x^{2} - 2 \, x + 8 $ & 6.2.4721393.4 & $ 2 $ & $ \infty $ & $\bigl\{ (1, 4), (1, 4), (1, 4), (3, 4) \bigr\}$\\
\lp & $ x^{6} - 2 \, x^{5} + x^{4} - 18 \, x^{3} + 18 \, x^{2} - 8 $ & 6.2.11279504.1 & $ 2 $ & $ \infty $ & $\bigl\{ (1, 4), (1, 4), (2, 1), (2, 1) \bigr\}$\\
\lp & $ x^{6} - x^{4} + 6 \, x^{2} - 8 $ & 6.2.5944352.3 & $ 2 $ & $ \infty $ & $\bigl\{ (1, 4), (1, 4), (2, 1), (2, 2) \bigr\}$\\
\lp & $ x^{6} - 3 \, x^{5} + 5 \, x^{4} - 5 \, x^{3} - 3 \, x^{2} + 5 \, x - 2 $ & 6.2.2150081.1 & $ 2 $ & $ \infty $ & $\bigl\{ (1, 4), (1, 4), (2, 2), (2, 2) \bigr\}$\\
\lp & $ x^{6} - x^{5} - 5 \, x^{4} + 5 \, x^{3} + x^{2} + 5 \, x - 4 $ & 6.2.946832.1 & $ 2 $ & $ \infty $ & $\bigl\{ (1, 4), (1, 4), (4, 1) \bigr\}$\\
\lp & $ x^{6} - 2 \, x^{5} - x^{4} + 2 \, x^{3} - x^{2} - 3 \, x + 2 $ & 6.2.365117.1 & $ 2 $ & $ \infty $ & $\bigl\{ (1, 4), (1, 4), (4, 2) \bigr\}$\\
\lp & $ x^{6} + 2 \, x^{5} - x^{4} - x^{3} - x^{2} - 2 $ & --- & $ 2 $ & $ \infty $ & $\bigl\{ (1, 4), (2, 1), (3, 4) \bigr\}$\\
\lp & $ x^{6} - x^{4} - 2 \, x^{2} - x - 1 $ & 6.2.255917.1 & $ 2 $ & $ \infty $ & $\bigl\{ (1, 4), (2, 2), (3, 4) \bigr\}$\\
\lp & $ x^{6} - 2 \, x^{5} + 3 \, x^{4} - x^{3} + 2 \, x - 1 $ & 6.2.89737.1 & $ 2 $ & $ \infty $ & $\bigl\{ (1, 4), (5, 4) \bigr\}$\\
\lp & $ x^{6} - 2 \, x^{5} - 5 \, x^{4} + 18 \, x^{3} - 10 \, x^{2} - 8 \, x + 8 $ & 6.2.48491968.1 & $ 2 $ & $ \infty $ & $\bigl\{ (2, 1), (2, 1), (2, 1) \bigr\}$\\
\lp & $ x^{6} - 2 \, x^{5} - x^{4} + 4 \, x^{3} + 2 \, x^{2} - 4 \, x - 8 $ & 6.2.2895824.1 & $ 2 $ & $ \infty $ & $\bigl\{ (2, 1), (2, 1), (2, 2) \bigr\}$\\
\lp & $ x^{6} - 2 \, x^{5} + 5 \, x^{4} - 4 \, x^{3} - 3 \, x^{2} + 6 \, x - 2 $ & 6.2.633788.1 & $ 2 $ & $ \infty $ & $\bigl\{ (2, 1), (2, 2), (2, 2) \bigr\}$\\
\lp & $ x^{6} - 4 \, x^{4} - 8 \, x^{3} - 4 \, x^{2} + 1 $ & 6.2.3548608.1 & $ 2 $ & $ \infty $ & $\bigl\{ (2, 1), (4, 1) \bigr\}$\\
\lp & $ x^{6} - x^{4} - 2 \, x^{3} + x^{2} + 2 \, x - 2 $ & 6.2.95852.1 & $ 2 $ & $ \infty $ & $\bigl\{ (2, 1), (4, 2) \bigr\}$\\
\lp & $ x^{6} - 2 \, x^{5} + 2 \, x^{3} - 2 \, x^{2} + 4 \, x + 1 $ & 6.2.332021.1 & $ 2 $ & $ \infty $ & $\bigl\{ (2, 2), (2, 2), (2, 2) \bigr\}$\\
\lp & $ x^{6} - 2 \, x^{5} + 2 \, x^{4} - 4 \, x^{3} + 2 \, x^{2} - 4 \, x + 1 $ & 6.2.242000.2 & $ 2 $ & $ \infty $ & $\bigl\{ (2, 2), (4, 1) \bigr\}$\\
\lp & $ x^{6} - x^{4} - x^{3} - x^{2} + 1 $ & 6.2.52441.1 & $ 2 $ & $ \infty $ & $\bigl\{ (2, 2), (4, 2) \bigr\}$\\
\lp & $ x^{6} - x^{5} - 3 \, x^{4} + 3 \, x^{3} + 3 \, x^{2} - x - 1 $ & 6.2.47081.1 & $ 2 $ & $ \infty $ & $\bigl\{ (3, 4), (3, 4) \bigr\}$\\
\lp & $ x^{6} + 2 \, x^{4} + x^{2} - 23 $ & 6.2.778688.2 & $ 2 $ & $ \infty $ & $\bigl\{ (6, 1) \bigr\}$\\
\lp & $ x^{6} - 2 \, x^{5} + 3 \, x^{3} - 2 \, x - 1 $ & 6.2.28037.1 & $ 2 $ & $ \infty $ & $\bigl\{ (6, 2) \bigr\}$\\
\lp\label{c6_GRH_50} & $ x^{6} + 8152 \, x^{5} + 131685 \, x^{4} + 3664 \, x^{3} + 49395 \, x^{2} + 496600 \, x + 23207 $ & --- & $ 4 $ & $ \infty $ & $\bigl\{ (1, 4), (1, 4), (1, 4), (1, 4), (1, 4), (1, 4) \bigr\}$\\
\lp & $ x^{6} + 474 \, x^{5} - 72 \, x^{4} - 424 \, x^{3} + 415 \, x^{2} - 434 \, x - 280 $ & --- & $ 4 $ & $ \infty $ & $\bigl\{ (1, 4), (1, 4), (1, 4), (1, 4), (2, 1) \bigr\}$\\
\lp\label{c6_GRH_52} & $ x^{6} + 120 \, x^{5} + 2 \, x^{4} - 169 \, x^{3} - 170 \, x^{2} + 176 \, x + 32 $ & --- & $ 4 $ & $ \infty $ & $\bigl\{ (1, 4), (1, 4), (1, 4), (1, 4), (2, 2) \bigr\}$\\
\lp & $ x^{6} + x^{5} + x^{4} - 21 \, x^{2} + 4 $ & --- & $ 4 $ & $ \infty $ & $\bigl\{ (1, 4), (1, 4), (1, 4), (3, 4) \bigr\}$\\
\lp\label{c6_GRH_54} & $ x^{6} - 484 \, x^{5} - 696 \, x^{4} + 346 \, x^{3} - 671 \, x^{2} - 538 \, x + 114 $ & --- & $ 4 $ & $ \infty $ & $\bigl\{ (1, 4), (1, 4), (2, 1), (2, 1) \bigr\}$\\
\lp & $ x^{6} - 15 \, x^{5} + 11 \, x^{4} + 25 \, x^{3} + 32 \, x^{2} + 34 \, x + 8 $ & --- & $ 4 $ & $ \infty $ & $\bigl\{ (1, 4), (1, 4), (2, 1), (2, 2) \bigr\}$\\
\lp & $ x^{6} - 4 \, x^{4} - 9 \, x^{2} + 4 $ & 6.4.13675204.1 & $ 4 $ & $ \infty $ & $\bigl\{ (1, 4), (1, 4), (2, 2), (2, 2) \bigr\}$\\
\lp & $ x^{6} - x^{5} - 7 \, x^{4} + 9 \, x^{3} + 11 \, x^{2} - 17 \, x + 2 $ & 6.4.3095536.1 & $ 4 $ & $ \infty $ & $\bigl\{ (1, 4), (1, 4), (4, 1) \bigr\}$\\
\lp & $ x^{6} - 5 \, x^{4} - 2 \, x^{3} + 7 \, x^{2} + 5 \, x - 2 $ & 6.4.1415907.1 & $ 4 $ & $ \infty $ & $\bigl\{ (1, 4), (1, 4), (4, 2) \bigr\}$\\
\lp & $ x^{6} + x^{5} - 3 \, x^{4} - 2 \, x^{3} - 2 \, x^{2} + x + 6 $ & --- & $ 4 $ & $ \infty $ & $\bigl\{ (1, 4), (2, 1), (3, 4) \bigr\}$\\
\lp & $ x^{6} - 3 \, x^{5} - 6 \, x^{4} + 12 \, x^{3} + 15 \, x^{2} + 3 $ & 6.4.32019867.1 & $ 4 $ & $ \infty $ & $\bigl\{ (1, 4), (2, 2), (3, 4) \bigr\}$\\
\lp & $ x^{6} - 2 \, x^{5} + x^{4} - 4 \, x^{2} + x + 2 $ & 6.4.321527.1 & $ 4 $ & $ \infty $ & $\bigl\{ (1, 4), (5, 4) \bigr\}$\\
\lp & $ x^{6} - 2 \, x^{5} - 10 \, x^{4} - 2 \, x^{3} - 11 \, x^{2} - 16 \, x + 16 $ & --- & $ 4 $ & $ \infty $ & $\bigl\{ (2, 1), (2, 1), (2, 1) \bigr\}$\\
\lp & $ x^{6} - 3 \, x^{5} - 2 \, x^{4} + 9 \, x^{3} - x^{2} - 4 \, x + 2 $ & 6.4.4293808.3 & $ 4 $ & $ \infty $ & $\bigl\{ (2, 1), (2, 1), (2, 2) \bigr\}$\\
\lp & $ x^{6} + 2 \, x^{4} - 6 \, x^{2} + 1 $ & 6.4.2149156.1 & $ 4 $ & $ \infty $ & $\bigl\{ (2, 1), (2, 2), (2, 2) \bigr\}$\\
\lp & $ x^{6} - 2 \, x^{5} - 5 \, x^{4} + 12 \, x^{3} - x^{2} - 6 \, x + 2 $ & 6.4.3182656.1 & $ 4 $ & $ \infty $ & $\bigl\{ (2, 1), (4, 1) \bigr\}$\\
\lp & $ x^{6} - x^{5} - 3 \, x^{4} + 5 \, x^{3} + 3 \, x^{2} - 4 \, x - 2 $ & 6.4.733588.1 & $ 4 $ & $ \infty $ & $\bigl\{ (2, 1), (4, 2) \bigr\}$\\
\lp & $ x^{6} - 5 \, x^{4} + 3 \, x^{2} + 2 $ & 6.4.1759688.1 & $ 4 $ & $ \infty $ & $\bigl\{ (2, 2), (2, 2), (2, 2) \bigr\}$\\
\lp & $ x^{6} - 3 \, x^{5} + x^{4} + 2 \, x^{3} - 3 \, x^{2} + 5 \, x + 1 $ & 6.4.478000.1 & $ 4 $ & $ \infty $ & $\bigl\{ (2, 2), (4, 1) \bigr\}$\\
\lp & $ x^{6} - 2 \, x^{5} - 2 \, x^{4} + 5 \, x^{3} - x^{2} - 3 \, x + 1 $ & 6.4.202375.1 & $ 4 $ & $ \infty $ & $\bigl\{ (2, 2), (4, 2) \bigr\}$\\
\lp & $ x^{6} - 3 \, x^{5} + x^{4} + 3 \, x^{3} - 3 \, x^{2} + x + 1 $ & 6.4.170471.1 & $ 4 $ & $ \infty $ & $\bigl\{ (3, 4), (3, 4) \bigr\}$\\
\lp & $ x^{6} - 3 \, x^{4} - 4 \, x^{3} + 2 \, x^{2} + 4 \, x + 1 $ & 6.4.526912.1 & $ 4 $ & $ \infty $ & $\bigl\{ (6, 1) \bigr\}$\\
\lp & $ x^{6} - x^{5} - 2 \, x^{4} + 3 \, x^{3} - x^{2} - 2 \, x + 1 $ & 6.4.92779.1 & $ 4 $ & $ \infty $ & $\bigl\{ (6, 2) \bigr\}$\\
\lp & $ x^{6} - 85 \, x^{4} - 16 \, x^{3} + 1156 \, x^{2} - 544 \, x + 64 $ & --- & $ 6 $ & $ \infty $ & $\bigl\{ (1, 4), (1, 4), (1, 4), (1, 4), (1, 4), (1, 4) \bigr\}$\\
\lp\label{c6_GRH_74} & $ x^{6} + 1028 \, x^{5} + 131059 \, x^{4} + 4194272 \, x^{3} + 8388639 \, x^{2} + 4194268 \, x + 131117 $ & --- & $ 6 $ & $ \infty $ & $\bigl\{ (1, 4), (1, 4), (1, 4), (1, 4), (2, 1) \bigr\}$\\
\lp & $ x^{6} + 496 \, x^{5} + 65611 \, x^{4} + 2097152 \, x^{3} + 2096311 \, x^{2} + 6032 \, x - 643 $ & --- & $ 6 $ & $ \infty $ & $\bigl\{ (1, 4), (1, 4), (1, 4), (1, 4), (2, 2) \bigr\}$\\
\lp & $ x^{6} + 8 \, x^{5} - 105 \, x^{4} + 8 \, x^{3} + 912 \, x^{2} + 16 \, x - 768 $ & --- & $ 6 $ & $ \infty $ & $\bigl\{ (1, 4), (1, 4), (1, 4), (3, 4) \bigr\}$\\
\lp & $ x^{6} + 2 \, x^{5} - 77 \, x^{4} + 4 \, x^{3} + 910 \, x^{2} + 4 \, x - 768 $ & --- & $ 6 $ & $ \infty $ & $\bigl\{ (1, 4), (1, 4), (2, 1), (2, 1) \bigr\}$\\
\lp\label{c6_GRH_78} & $ x^{6} + 4 \, x^{5} - 105 \, x^{4} + 4 \, x^{3} + 910 \, x^{2} + 4 \, x - 768 $ & --- & $ 6 $ & $ \infty $ & $\bigl\{ (1, 4), (1, 4), (2, 1), (2, 2) \bigr\}$\\
\lp & $ x^{6} - 3 \, x^{5} - 9 \, x^{4} + 23 \, x^{3} + 20 \, x^{2} - 32 \, x - 8 $ & 6.6.56269193.1 & $ 6 $ & $ \infty $ & $\bigl\{ (1, 4), (1, 4), (2, 2), (2, 2) \bigr\}$\\
\lp & $ x^{6} - 9 \, x^{4} + 13 \, x^{2} - 4 $ & 6.6.24167056.1 & $ 6 $ & $ \infty $ & $\bigl\{ (1, 4), (1, 4), (4, 1) \bigr\}$\\
\lp & $ x^{6} - 2 \, x^{5} - 7 \, x^{4} + 14 \, x^{3} + 5 \, x^{2} - 13 \, x + 4 $ & 6.6.10330877.1 & $ 6 $ & $ \infty $ & $\bigl\{ (1, 4), (1, 4), (4, 2) \bigr\}$\\
\lp & $ x^{6} - 2 \, x^{5} - 6 \, x^{4} + 11 \, x^{3} + 7 \, x^{2} - 11 \, x - 2 $ & 6.6.20413244.1 & $ 6 $ & $ \infty $ & $\bigl\{ (1, 4), (2, 1), (3, 4) \bigr\}$\\
\lp & $ x^{6} - x^{5} - 7 \, x^{4} + 4 \, x^{3} + 13 \, x^{2} - 2 \, x - 4 $ & 6.6.7432373.1 & $ 6 $ & $ \infty $ & $\bigl\{ (1, 4), (2, 2), (3, 4) \bigr\}$\\
\lp & $ x^{6} - 6 \, x^{4} - x^{3} + 8 \, x^{2} + x - 2 $ & 6.6.1868969.1 & $ 6 $ & $ \infty $ & $\bigl\{ (1, 4), (5, 4) \bigr\}$\\
\lp & $ x^{6} - 2 \, x^{5} - 17 \, x^{4} + 54 \, x^{3} - 34 \, x^{2} - 8 \, x + 8 $ & 6.6.31554496.1 & $ 6 $ & $ \infty $ & $\bigl\{ (2, 1), (2, 1), (2, 1) \bigr\}$\\
\lp & $ x^{6} - 14 \, x^{4} - 2 \, x^{3} + 15 \, x^{2} - 2 \, x - 2 $ & 6.6.34350464.1 & $ 6 $ & $ \infty $ & $\bigl\{ (2, 1), (2, 1), (2, 2) \bigr\}$\\
\lp & $ x^{6} - 3 \, x^{5} - 6 \, x^{4} + 13 \, x^{3} + 9 \, x^{2} - 12 \, x - 4 $ & 6.6.35478972.1 & $ 6 $ & $ \infty $ & $\bigl\{ (2, 1), (2, 2), (2, 2) \bigr\}$\\
\lp & $ x^{6} - 9 \, x^{4} - 6 \, x^{3} + 15 \, x^{2} + 14 \, x + 2 $ & 6.6.9186752.1 & $ 6 $ & $ \infty $ & $\bigl\{ (2, 1), (4, 1) \bigr\}$\\
\lp & $ x^{6} - x^{5} - 6 \, x^{4} + 4 \, x^{3} + 7 \, x^{2} - 2 \, x - 2 $ & 6.6.3072812.1 & $ 6 $ & $ \infty $ & $\bigl\{ (2, 1), (4, 2) \bigr\}$\\
\lp & $ x^{6} - 12 \, x^{4} - 4 \, x^{3} + 16 \, x^{2} + 4 \, x - 1 $ & 6.6.6555125.1 & $ 6 $ & $ \infty $ & $\bigl\{ (2, 2), (2, 2), (2, 2) \bigr\}$\\
\lp & $ x^{6} - 7 \, x^{4} + 12 \, x^{2} - 2 $ & 6.6.3195392.1 & $ 6 $ & $ \infty $ & $\bigl\{ (2, 2), (4, 1) \bigr\}$\\
\lp & $ x^{6} - x^{5} - 7 \, x^{4} + 7 \, x^{3} + 12 \, x^{2} - 12 \, x - 1 $ & 6.6.1312625.1 & $ 6 $ & $ \infty $ & $\bigl\{ (2, 2), (4, 2) \bigr\}$\\
\lp & $ x^{6} - x^{5} - 7 \, x^{4} + 9 \, x^{3} + 7 \, x^{2} - 9 \, x - 1 $ & 6.6.905177.1 & $ 6 $ & $ \infty $ & $\bigl\{ (3, 4), (3, 4) \bigr\}$\\
\lp & $ x^{6} - 7 \, x^{4} + 14 \, x^{2} - 7 $ & 6.6.1075648.1 & $ 6 $ & $ \infty $ & $\bigl\{ (6, 1) \bigr\}$\\
\lp & $ x^{6} - x^{5} - 7 \, x^{4} + 2 \, x^{3} + 7 \, x^{2} - 2 \, x - 1 $ & 6.6.300125.1 & $ 6 $ & $ \infty $ & $\bigl\{ (6, 2) \bigr\}$\\
\hline
\end{longtable}
\end{landscape}

\end{document}